\documentclass[12pt]{article}
\usepackage{mathrsfs}
\usepackage{CJK}
\usepackage{bbm}
\usepackage{stmaryrd}
\usepackage{shortvrb,psfrag}
\usepackage{enumerate}
\usepackage{amsmath}
\usepackage{cases}
\usepackage{amsfonts}
\usepackage{latexsym,amssymb,amsmath,mathrsfs,amsbsy, amsthm}
\usepackage[dvips]{graphicx}
\usepackage[usenames]{color}
\usepackage{xspace,colortbl}
\usepackage{pdfsync}
\usepackage[dvipdfm,
pdfstartview=FitH, CJKbookmarks=true, bookmarksnumbered=true,
bookmarksopen=true,
citecolor=blue
,colorlinks=true
,pdfborder=1
,pdfstartpage=2
,pdfstartview=Fit]{hyperref}
\allowdisplaybreaks

\newtheorem{thm}{Theorem}[section]
\newtheorem{lem}[thm]{Lemma}

\newtheorem{rem}{Remark}
\newtheorem{cor}[thm]{Corollary}

\newcommand{\beq}{\begin{equation}}
\newcommand{\eeq}{\end{equation}}
\newcommand{\ben}{\begin{eqnarray}}
\newcommand{\een}{\end{eqnarray}}
\newcommand{\beno}{\begin{eqnarray*}}
\newcommand{\eeno}{\end{eqnarray*}}

\begin{document}
\title{{\bf Stability of the Couette flow under the 2D steady Navier-Stokes flow
}
\\Wendong WANG
\\[2mm]
{\small  School of Mathematical Sciences, Dalian University of Technology, Dalian 116024,  China}\\
{\small Email: wendong@dlut.edu.cn}\\[2mm]}

\date{\today}
\maketitle


\begin{abstract}
In this note, we investigate the stability property of shear flows under the 2D stationary Navier-Stokes equations, and we obtain that the Couette flow $(y,0)$ is stable under the space of $\mathcal{D}^{1,q}(\mathbb{R}^2)$ for any $1<q<\infty$  and unstable in the space of $\mathcal{D}^{1,\infty}(\mathbb{R}^2)$. A key observation is the anisotropic cut-off function. We also consider the Poiseuille flow $(y^2,0)$, which is stable in  $\mathcal{D}^{1,q}(\mathbb{R}^2)$ with $\frac43<q\leq4.$
\end{abstract}


\setcounter{equation}{0}
\section{Introduction}

Consider the incompressible steady Navier-Stokes equations in a domain $\Omega\subset \mathbb{R}^2$:
\begin{equation}\label{eq:NS-2D}
\left\{\begin{array}{llll}
-\mu\Delta u+u\cdot \nabla u+\nabla \pi=0,\\
{\rm div }~ u=0,
\end{array}\right.
\end{equation}
where $\mu$ denotes the viscosity coefficient. We assume $\mu=1$ for simplicity.

One fundamental question is to investigate the well-posedness property of (\ref{eq:NS-2D}). The existence on an exterior domain attracts the attention of many mathematicians  when the boundary condition is given at infinity:
\ben\label{eq:boundary infty}
\lim_{|x|\rightarrow\infty}u(x)=u_\infty,
\een
where $u_\infty$ is a constant vector, for example, see Leray \cite{Leray} and Russo \cite{Ru2009}.
They constructed a solution whose Dirichlet energy is bounded:
\ben\label{eq:energy bound-NS-2D}
D(u)=\int_{R^2}|\nabla u|^2dx<\infty
\een
but it's difficulty to verify that it satisfies the condition (\ref{eq:boundary infty}). Hence, one challenging problem is to prove the constructed solution satisfying the asymptotic behavior at $\infty$.
Gilbarg-Weinberger \cite{GW1978} described the asymptotic behavior of the velocity, the pressure and the vorticity, where they showed that $u(x)=o(\ln|x|)$ and
\beno
\lim_{r\rightarrow\infty}\int_0^{2\pi}|u(r,\theta)-\bar{u}|^2d\theta=0
\eeno
for some constant vector $\bar{u}$.
Later, Amick \cite{Amick} proved that $u\in L^\infty$ under zero
boundary condition.  Recently, Korobkov-Pileckas-Russo in \cite{KPR2017} and \cite{KPR2018} obtained that
\beno
\lim_{|x|\rightarrow\infty}u(x)=\bar{u}.
\eeno
More references on the existence and asymptotic behavior of solutions in an exterior domain, we refer to \cite{GNP1997,Ru2010,GG2011,PR2012,KPR2014,DI2017} and the references therein.

When $\Omega$ is the whole space, an interesting question is to study the classification of solutions of  (\ref{eq:NS-2D}). In details,  we are concerned on the solution spaces of (\ref{eq:NS-2D}), or Liouville properties around some special solutions such as shear flows.
The shear flow is like the form of $U(y)=(g(y), 0)$, and it follows from (\ref{eq:NS-2D}) that $g(y)=c, y,$ or $y^2$. As in \cite{Galdi}, for a domain $\Omega$ and $1\leq q\leq \infty$,  we define the following linear space (without
topology)
\beno
\mathcal{D}^{1,q}(\Omega)\doteq\{u\in L_{loc}^1(\Omega), \nabla u\in L^q(\Omega)\},
\eeno
which describes the growth of the energy.
Furthermore, for $\alpha\in [0,1]$ and $1\leq p\leq \infty$, we introduce another space
\beno
\chi^{\alpha,p}(\Omega)\doteq\{u\in L_{loc}^1(\Omega), \frac{|u|}{(1+|x|)^\alpha}\in L^{p}(\Omega)\},
\eeno
which describes the growth of $u$. Obviously, $\chi^{0,p}(\Omega)$ is the usual $L^p(\Omega)$ space.

For $g(y)=c$, let us recall some known results on
 this issue.   Under the condition
 (\ref{eq:boundary infty}), the smooth solution $u$ is indeed bounded and a Liouville theorem being more in the spirit of the classical one for entire analytic functions
was obtained by Koch-Nadirashvili-Seregin-Sverak \cite{KNSS} as a byproduct of their
work on the nonstationary case. If $u\in \chi^{0,p}(\mathbb{R}^2)$ for $p>1$, then $u$ is trivial, see Zhang \cite{ZG2015}. As suggested by Fuchs-Zhong in \cite{FZ2011}, the stable space may be $\chi^{\alpha,\infty}(\mathbb{R}^2)$ with $\alpha<1$ as the property of harmonic functions, since the linear solutions are the counterexamples; see also Yau \cite{Yau} and Peter Li-Tam \cite{Li-Tam}, where they considered the space of harmonic functions on complete manifold with
nonnegative Ricci curvature with linear growth. When $\alpha\in [0, 1/7)$ and $u\in \chi^{\alpha,\infty}(\mathbb{R}^2)$, $u$ is a constant vector by Fuchs-Zhong \cite{FZ2011}. The component is improved to $\alpha<\frac13$ with help of the vorticity equation  by Bildhauer-Fuchs-Zhang in \cite{BFZ2013}.

On the other hand, for the growth of $\nabla u$, Gilbarg-Weinberger proved the above Liouville type theorem by assuming (\ref{eq:energy bound-NS-2D}) in  \cite{GW1978}, where they made use of the fact that the vorticity function satisfies
a nice elliptic equation to which a maximum principle applies.
The assumption on boundedness of the Dirichlet energy can be relaxed to $\nabla u\in L^p(\mathbb{R}^2)$ with some $p\in (\frac65,3]$, see Bildhauer-Fuchs-Zhang \cite{BFZ2013} for generalized Navier-Stokes equations. If $u\in \mathcal{D}^{1,q}(\mathbb{R}^2)$ for $1<p<\infty$, the constant $u$ follows by the author in \cite{Wang}.
The above results also can be generalized to the shear thickening flows, for example see \cite{Fu2012exist,Fu2012Liou,FZ2012,ZG2013,JK2014}.
For the two dimensional steady MHD equations, the similar Liouvile type theorems was obtained by Y. Wang and the author in \cite{WW} by assuming (\ref{eq:energy bound-NS-2D}) or $u\in \chi^{0,p}(\mathbb{R}^2)$ with $2<p\leq \infty$, where the smallness conditions of the magnetic field are added. See also the recent result  in \cite{Wang} for $u\in \chi^{\alpha,\infty}(\mathbb{R}^2)$ with $\alpha<\frac13$ by using the idea of \cite{BFZ2013} and energy estimates in an annular domain.

Next we consider the stable space of shear flows in $\mathcal{D}^{1,q}$ or $\chi^{\alpha,\infty}.$
Let $(u,\pi)$ be a smooth solution of (\ref{eq:NS-2D})
and  the vorticity $\tilde{w}=\partial_2u_1-\partial_1 u_2$, then the vorticity equations are as follows:
\begin{equation}\label{eq:NSv-2D}
-\Delta \tilde{w}+u\cdot \nabla \tilde{w}=0,
\end{equation}

We will first study the stability of the Couette flow $U = (y, 0)$, which is a solution of (\ref{eq:NS-2D}). Let $v = u-U$ be the perturbation of the
velocity satisfying
\begin{equation}\label{eq:NS-2D-c}
\left\{\begin{array}{llll}
-\Delta v+v\cdot \nabla v+\nabla \pi+(v_2,0)
+y\partial_x v=0,\\
{\rm div }~ v=0,
\end{array}\right.
\end{equation}
Let $w=\partial_2v_1-\partial_1 v_2$, then
\begin{equation}\label{eq:NSv-2D-c}
-\Delta w+v\cdot \nabla w+y\partial_xw=0.
\end{equation}

Now we state our main result on the Couette flow:
\begin{thm}\label{thm:generalization of GW}
Let $(u,\pi)$ be a smooth solution of the 2D Navier-Stokes equations (\ref{eq:NS-2D}) defined over the entire plane. For $U = (y, 0)$, assume that $v=u-U\in \mathcal{D}^{1,q}(\mathbb{R}^2)$ with $1< q< \infty$. Then $v$ and $\pi$ are constants.
\end{thm}

\begin{rem}
Obviously, the above result fails in the space of $\mathcal{D}^{1,\infty}(\mathbb{R}^2)$, since the linear solutions are not unique. The above result also shows that the stable spaces is similar as the constant solution(see \cite{Wang}).
It is worth mentioning that the stability threshold in Sobolev spaces for the 2D time-dependent Navier-Stokes is more complicated, for example, see Bedrossian-Germain-Masmoudi \cite{BGM}, Bedrossian-Wang-Vicol
\cite{BVF} and Chen-Li-Wei-Zhang \cite{CLWZ}, where if the initial velocity is around the Couette
flow
\beno
\|u_0-(y,0)\|_{H^2}\leq c Re^{-\frac12}
\eeno
for a small $c$, then the solution still stays in this space for any time.
\end{rem}


If the velocity is largely growing around the Couette flow, we have the following stability estimate:
\begin{thm}\label{thm:generalization of KNSS-C}
Let $(u,\pi)$ be a smooth solution of the 2D Navier-Stokes equations (\ref{eq:NS-2D}) defined over the entire plane and $v=u-U$ satisfies the growth estimates $v\in \chi^{\alpha,\infty}(\mathbb{R}^2)$ for $0<\alpha<\frac15$, where $U = (y, 0)$. Then $v$ and $\pi$ are constants.
\end{thm}

\begin{rem}
The above result generalized the Liouville type theorems around the trivial solution in
\cite{KNSS,BFZ2013} to the Couette flow.
\end{rem}

The similar arguments can applied to the Poiseuille flow, which is stated as follows.
\begin{cor}\label{thm:generalization of KNSS-P}
Let $(u,\pi)$ be a smooth solution of the 2D Navier-Stokes equations (\ref{eq:NS-2D}) defined over the entire plane. For $U = (y^2, 0)$, let $v=u-U$. Then $v$ and $\pi$ are constants, if one of the following conditions holds:\\
(i)  $v\in \mathcal{D}^{1,q}(\mathbb{R}^2)$ with $\frac43< q\leq 4$;\\
(ii) $v$ satisfies the growth estimates $v\in \chi^{\alpha,\infty}(\mathbb{R}^2)$ for $0<\alpha<\frac18$.
\end{cor}


Let us recall a result of Gilbarg-Weinberger in \cite{GW1978} about the decay of functions with finite Dirichlet integrals.
\begin{lem}[Lemma 2.1, 2.2, \cite{GW1978}]
\label{lem:GW}
Let a $C^1$ vector-valued function $f(x)=(f_1,f_2)(x)=f(r,\theta)$ with $r=|x|$ and $x_1=r\cos\theta$.  There holds finite Dirichlet integral in the range $r>r_0$, that is
\beno
\int_{r>r_0}|\nabla f|^2\,dxdy<\infty.
\eeno
Then, we have
\beno
\lim_{r\rightarrow\infty} \frac{1}{\ln r}\int_0^{2\pi}|f(r,\theta)|^2d\theta=0.
\eeno
\end{lem}

If, furthermore, we assume $\nabla f \in L^p(\mathbb{R}^2)$ for some $2< p < \infty$,
then the above decay property can be improved to be point-wise uniformly.
More precisely, we have

\begin{lem}[Theorem II.9.1 \cite{Galdi}]
\label{lem:Galdi}
Let $\Omega \subset \mathbb{R}^2$ be an exterior domain. \\
(i) Let
\[
\nabla f \in L^2 \cap L^p (\Omega),
\]
for some $2< p < \infty$.
Then
\[
\lim_{|x| \to \infty} \frac{|f(x)|}{\sqrt{ \ln (|x|)}} = 0,
\]
uniformly.\\
(ii)Let
\beno
\nabla f \in  L^p (\Omega),
\eeno
for some $2< p < \infty$.
Then
\beno
\lim_{|x| \rightarrow \infty} \frac{|f(x)|}{|x|^{\frac{p-2}{p}}} = 0,
\eeno
uniformly.
\end{lem}

%

Throughout this article, $C(\alpha_1,\cdots,\alpha_n)$ denotes a constant depending on $\alpha_1,\cdots,\alpha_n$, which may be different from line to line.

\section{Proof of Theorem \ref{thm:generalization of GW}}

In this section, we are aimed to prove Theorem \ref{thm:generalization of GW}. First, let us  prove a similar result as Gilbarg-Weinberger in \cite{GW1978} about the decay of functions with finite $\mathcal{D}^{1,q}$ gradient integrals.
\begin{lem}
\label{lem:GW2}
Let a $C^1$ vector-valued function $f(x)=(f_1,f_2)(x)=f(r,\theta)$ with $r=|x|$ and $x_1=r\cos\theta$.  There holds
\beno
\int_{r>r_0}|\nabla f|^q\,dx<\infty,\quad 1<q<2.
\eeno
Then, we have
\beno
\limsup_{r\rightarrow\infty} \int_0^{2\pi}|f(r,\theta)|^qd\theta<\infty.
\eeno
\end{lem}

{\bf Proof of  Lemma \ref{lem:GW2}.}
By H\"{o}lder inequality we have
\beno
\frac{d}{dr}\left(\int_0^{2\pi}|f(r,\theta)|^qd\theta\right)^{\frac1q}&\leq&\left(\int_0^{2\pi}|f|^qd\theta\right)^{\frac1q-1}\int_0^{2\pi}|f|^{q-1}|f_r|d\theta\\
&\leq&\left(\int_0^{2\pi}|f_r|^{q}d\theta\right)^{\frac1q}
\eeno
Integrating from $r_1$ with $r_1\geq r_0$, we get
\beno
&&\left(\int_0^{2\pi}|f(r,\theta)|^qd\theta\right)^{\frac1q}-\left(\int_0^{2\pi}|f(r_1,\theta)|^qd\theta\right)^{\frac1q}\\
&\leq&
\int_{r_1}^r\left(\int_0^{2\pi}|f_r|^{q}d\theta\right)^{\frac1q} dr\\
&\leq& \left(\int_{r>r_0}|\nabla f|^q\,dx\right)^{\frac1q}\left(\int_{r_1}^rr^{-\frac{1}{q-1}}dr\right)^{1-\frac1q}\\
&\leq& \left(\int_{r>r_0}|\nabla f|^q\,dx\right)^{\frac1q}\left(\frac{q-1}{2-q}\right)^{1-\frac1q}r_1^{\frac{q-2}{q}}
\eeno
which yields the required result.

{\bf Proof of  Theorem \ref{thm:generalization of GW}.}

{\bf Step I. Case of $2<q<\infty.$}
Let $\eta(x,y)\in C_0^\infty(\mathbb{R}^2)$ be a cut-off function with $0\leq \eta\leq 1$ satisfying $\eta(x,y)=\eta_1(x)\eta_2(y)$, where
\begin{align*} \eta_1(x)=\left\{
\begin{aligned}
&1,\quad |x|\leq R^{\beta},\\
&0, \quad |x|>2R^{\beta}.
\end{aligned}
\right. \end{align*}
where $1<\beta<(1-\frac2q)^{-1}$,
and
\begin{align*} \eta_2(y)=\left\{
\begin{aligned}
&1,\quad |y|\leq R,\\
&0, \quad |y|>2R.
\end{aligned}
\right. \end{align*}

Multiply $q\eta|w|^{q-2}w$  on both sides of (\ref{eq:NSv-2D-c}), and we have
\ben\label{eq:energy estimate-w}
I&\doteq&\frac{4(q-1)}{q}\int_{\mathbb{R}^2}|\nabla (|w|^{\frac{q}{2}})|^2 \eta dxdy\nonumber\\
&\leq&\int_{\mathbb{R}^2} |w|^{q}\triangle \eta dxdy+\int_{\mathbb{R}^2} |w|^{q}y\partial_x\eta dxdy\nonumber\\
&&+\int_{\mathbb{R}^2} |w|^{q}v\cdot \nabla \eta dxdy\doteq I_1+I_2+I_3
\een
Since $ v\in \mathcal{D}^{1,q}$ and $\beta>1$,
obviously $I_1\rightarrow 0$  and
\beno
I_2\leq C R^{1-\beta}\rightarrow 0,
\eeno
as $R\rightarrow\infty.$
About the term $I_3$, due to Lemma \ref{lem:Galdi}, for large $R>0$ we have
\beno
|v(x,y)|\leq |(x,y)|^{1-\frac{2}{q}}
\eeno
Thus we have
\beno
I_3\leq C R^{\beta(1-\frac1q)-1}\rightarrow 0,
\eeno
as $R\rightarrow\infty,$ since
\beno
\beta(1-\frac2q)<1.
\eeno

Consequently, we get
$
\nabla (|w|^{\frac{q}{2}})\equiv0,
$
which implies that
$
w\equiv C.
$
Due to ${\rm div}~ v=0$, it follows that
\beno
\triangle v\equiv 0
\eeno
which and the known condition $ v\in \mathcal{D}^{1,q}$ yield that
\beno
\nabla v\equiv 0.
\eeno
Hence $v$ and $\pi$ are constant.

{\bf Step II. Case of $1<q\leq 2.$}
We take a cut-off function $\phi$ as follows.
\begin{enumerate}
\item[i).] Let $r=\sqrt{x^2+y^2}$.
$\phi$ is radially decreasing and satisfies
\begin{align*} \phi(x,y)=\phi(r)=\left\{
\begin{aligned}
&1,\quad r\leq \rho,\\
&0, \quad r\geq\tau,
\end{aligned}
\right. \end{align*}
where $0<\frac{R}{2}\leq\frac23\tau<\frac34 R\leq  \rho<\tau\leq R$;

\item[ii).]
 $|\nabla\phi|(x,y)\leq \frac{C}{\tau-\rho}$ for all $(x,y)\in \mathbb{R}^2$.

\end{enumerate}

 Multiplying both sides of
\eqref{eq:NSv-2D-c} by $\phi w$ respectively and then
 applying integration by parts, we arrive at
\begin{align}
\label{eq:vorticity equ}
& \hspace{-6mm} \int_{\mathbb{R}^2} \phi |\nabla w |^2 \,dxdy \notag\\
& = -  \int_{\mathbb{R}^2} \nabla w  \cdot \nabla \phi w  \,dxdy
+ \frac12 \int_{\mathbb{R}^2} v \cdot \nabla  \phi w^2 \,dxdy + \frac12 \int_{\mathbb{R}^2}  y\partial_x  \phi w^2 \,dxdy \notag\\
&\doteq I_1' +I_2'+ I_3'.
\end{align}
In what follows we shall estimate $I_j'$ for $j = 1,2,3$ one by one.

For the term $I_1'$, by H\"older's inequality  we have
\begin{align*}
I_1' \le \frac{C}{\tau-\rho} \| \nabla w \|_{L^2 (B_\tau)} \| w \|_{L^2 (B_\tau )}
\end{align*}
Using the following Poincar\'{e}-Sobolev inequality(see, for example, Theorem 8.11 and 8.12 \cite{LL})
\ben
\label{eq:poincare-sobolev}
\|w\|_{L^2(B_\tau)}\leq C \|\nabla w\|_{L^2(B_\tau)}^{1-\frac q2}\|w\|_{L^q(B_\tau)}^{\frac q2}+C\tau^{1-\frac2q}\|w\|_{L^q(B_\tau)},
\een
which yields that
\ben\label{eq:estimate of I1}
I_1'
\leq \frac18 \int_{B_{\tau}} |\nabla w|^2 \,dx + \frac {C}{(\tau-\rho)^{\frac4q}}+\frac {C\tau^{2-\frac4q}}{(\tau-\rho)^2},
\een
by noting that $\|w\|_{L^q(B_\tau)}\leq \|\nabla v\|_{L^q(\mathbb{R}^2)}<\infty.$

For the terms $I_2'$,
let
\beno
\bar{f}(r)=\frac{1}{2\pi}\int_0^{2\pi}f(r,\theta)d\theta,
\eeno
then by Wirtinger's inequality (for example, for $p=2$ see Chapter II.5 \cite{Galdi}) we have
\ben\label{eq:Wirtinger}
\int_0^{2\pi}|f-\bar{f}|^p \, d\theta\leq C(p)\int_0^{2\pi}|\partial_\theta f|^pd\theta,
\een
for $1\leq p<\infty$.

Then by using (\ref{eq:Wirtinger}), Lemma \ref{lem:GW} and Lemma \ref{lem:GW2} we have
\begin{align*}
I_2'
&\leq \left|\int_{\mathbb{R}^2}w^2\,(v-\bar{v})\cdot\nabla\phi \,dxdy\right|+\left|\int_{\mathbb{R}^2}w^2\,\bar{v}\cdot\nabla\phi \,dxdy\right|\\
&   \leq \frac{ C}{\tau-\rho}\left( \int_{B_\tau}w^{2q'}\right)^{\frac1{q'}}
\left( \int_{\frac{R}2<r<R} \int_0^{2\pi}|v(r,\theta)- \bar{v} |^q \, d\theta \,r dr \right)^{\frac1q}\\
&\hspace{8mm} + \frac{C}{\tau-\rho}  \int_{B_\tau\setminus B_{\frac{R}2}} w^2 \left( \int_0^{2\pi} |v (r,\theta)|^q \, d\theta \right)^{\frac1q} \,dxdy\\
&   \leq \frac{ CR}{\tau-\rho}\left( \int_{B_\tau}w^{2q'}\right)^{\frac1{q'}}\left( \int_{\frac{R}2<r<R}\frac{1}{r^q}  \int_0^{2\pi}|\partial_\theta v|^qd\theta \,rdr\right)^{\frac1q}\\
&\hspace{8mm} +C \frac{(\ln R)^{\frac12}}{\tau-\rho} \int_{B_\tau}w^2\,dxdy.
\end{align*}
Using Poincar\'{e}-Sobolev inequality again,
\ben
\label{eq:poincare-sobolev2}
\|w\|_{L^{2q'}(B_\tau)}\leq C \|\nabla w\|_{L^2(B_\tau)}^{1-\frac{q}{2q'}}\|w\|_{L^q(B_\tau)}^{\frac{q}{2q'}}+C\tau^{1-\frac3q}\|w\|_{L^q(B_\tau)},
\een
which and (\ref{eq:poincare-sobolev}) imply that
\ben\label{eq:estimate of I2}
I_2'& \leq &\frac18\left( \int_{B_\tau}|\nabla w|^2\right)+C\left(\frac{ R}{\tau-\rho}\|\nabla v\|_{L^q(B_R\backslash B_{R/2})}\right)^{\frac{2q'}{q}}+\frac{ CR\tau^{2-\frac6q}}{\tau-\rho}\nonumber\\
&&+C\left( \frac{\sqrt{\ln R}}{\tau-\rho}\right)^{\frac2q}+C\left( \frac{\sqrt{\ln R}}{\tau-\rho}\right)\tau^{2-\frac4q} ,
\een
where we used the boundedness of $\mathcal{D}^{1,q}$ integral.

For the term $I_3'$, we have
\beno
I_3'\leq\left|\int_{\mathbb{R}^2}
 w^2 \,y\partial_x\phi \,dxdy\right|  \leq \frac{ CR}{\tau-\rho}\left( \int_{B_\tau\backslash{B_{\frac23\tau}}}w^{2}dxdy\right)
\eeno
and using Poincar\'{e}-Sobolev inequality in a cylinder domain, a slightly different version of (\ref{eq:poincare-sobolev}) is
\beno
\label{eq:poincare-sobolev3}
\|w\|_{L^2(B_\tau\backslash{B_{\frac23\tau}})}\leq C \|\nabla w\|_{L^2(B_\tau\backslash{B_{\frac23\tau}})}^{1-\frac q2}\|w\|_{L^q(B_\tau\backslash{B_{\frac23\tau}})}^{\frac q2}+C\tau^{1-\frac2q}\|w\|_{L^q(B_\tau\backslash{B_{\frac23\tau}})},
\eeno
which implies that
\ben\label{eq:estimate of I3}
I_3'
&\leq & \frac18\left( \int_{B_\tau}|\nabla w|^2\right)+C \left(\frac{R}{\tau-\rho}\right)^{\frac2q}\|w\|_{L^q(B_\tau\backslash{B_{\frac23\tau}})}^2\nonumber\\
&&+ C \frac{R}{\tau-\rho}\tau^{2-\frac4q}\|w\|_{L^q(B_\tau\backslash{B_{\frac23\tau}})}^2.
\een

Collecting the estimates of $I_1,\cdots, I_3$, by (\ref{eq:estimate of I1}), (\ref{eq:estimate of I2}) and (\ref{eq:estimate of I3}) we have
\beno
 &&\int_{B_\rho}|\nabla w|^2dxdy\\&\leq& \frac12 \int_{B_\tau}|\nabla w|^2+ \frac {C}{(\tau-\rho)^{\frac4q}}+\frac {C\tau^{2-\frac4q}}{(\tau-\rho)^2}++\frac{ CR\tau^{2-\frac6q}}{\tau-\rho}\\
&&+C\left( \frac{\sqrt{\ln R}}{\tau-\rho}\right)^{\frac2q}+C\left( \frac{\sqrt{\ln R}}{\tau-\rho}\right)\tau^{2-\frac4q}+C\left(\frac{ R}{\tau-\rho}\|\nabla v\|_{L^q(B_R\backslash B_{R/2})}\right)^{\frac{2q'}{q}}\\  &&+C \left(\frac{R}{\tau-\rho}\right)^{\frac2q}\|w\|_{L^q(B_R\backslash{B_{R/2}})}^2+ C \frac{R}{\tau-\rho}\tau^{2-\frac4q}\|w\|_{L^q(B_R\backslash{B_{R/2}})}^2
\eeno
Then an application of Giaquinta's iteration lemma \cite[Lemma 3.1]{G83} yields
\beno
 &&\int_{B_{R/2}}|\nabla w|^2dxdy\leq CR^{-\frac{4}{q}}+C\left( \frac{\sqrt{\ln R}}{R}\right)\\
&&+C\left(\|\nabla v\|_{L^q(B_R\backslash B_{R/2})}\right)^{\frac{2q'}{q}}+C(R^{2-\frac4q}+1
  ) \|w\|_{L^q(B_R\backslash{B_{R/2}})}^2
\eeno
Letting $R\rightarrow\infty$, we have
\beno
\nabla w\equiv 0
\eeno
and $w\equiv C$. Similar arguments as in Step I, we complete the proof.

\section{Proof of Theorem \ref{thm:generalization of KNSS-C}}

In this section, we will prove Theorem \ref{thm:generalization of KNSS-C}.

{\bf Proof.}
Let $\eta(x,y)\in C_0^\infty(\mathbb{R}^2)$ be a cut-off function on a cylinder domain with $0\leq \eta\leq 1$ satisfying $\eta(x,y)=\eta_1(x)\eta_2(y)$, where
\begin{align*} \eta_1(x)=\left\{
\begin{aligned}
&1,\quad |x|\leq R^{\beta},\\
&0, \quad |x|>2R^{\beta}.
\end{aligned}
\right. \end{align*}
where $\beta> 1$, to be decided,
and
\begin{align*} \eta_2(y)=\left\{
\begin{aligned}
&1,\quad |y|\leq R,\\
&0, \quad |y|>2R.
\end{aligned}
\right. \end{align*}

Write $w^{2q}=(w^2)^{q}$. As in \cite{BFZ2013}(see also \cite{Wang}), for $q\geq 2,\ell\geq q$, we have
\beno\label{eq:w estimate-mhd}
\int_{\mathbb{R}^2}w^{2q}\eta^{2\ell}dxdy\nonumber&=& \int_{\mathbb{R}^2}(\partial_2 v_1-\partial_1 v_2)w^{2q-2}w\eta^{2\ell}dxdy\nonumber\\
&=&  \int_{\mathbb{R}^2}(v_2,-v_1)\cdot \nabla [w^{2q-2}w\eta^{2\ell}]dxdy\nonumber\\
&\leq &(2q-1)\int_{\mathbb{R}^2}|v||\nabla w| w^{2q-2}\eta^{2\ell}dxdy+2\ell\int_{\mathbb{R}^2}|v|| \nabla \eta| |w|^{2q-1}\eta^{2\ell-1}dxdy\nonumber\\
&\leq&\frac12 \int_{\mathbb{R}^2}w^{2q}\eta^{2\ell}dxdy+C(q)\int_{\mathbb{R}^2}|v|^2|\nabla w|^2 w^{2q-4}\eta^{2\ell}dxdy\nonumber\\
&&+2\ell\int_{\mathbb{R}^2}|v|| \nabla \eta| |w|^{2q-1}\eta^{2\ell-1}dxdy
\eeno

Due to the growth estimates $v\in \chi^{\alpha,\infty}$,  we have
\ben\label{eq:w estimate-mhd4}
\int_{\mathbb{R}^2}w^{2q}\eta^{2\ell}dxdy
&\leq& C(q)R^{2\alpha\beta}\int_{\mathbb{R}^2}|\nabla w|^2 w^{2q-4}\eta^{2\ell}dxdy\nonumber\\
&&+ C(\ell)R^{\alpha\beta-1}\int_{\mathbb{R}^2} |w|^{2q-1}\eta^{2\ell-1}dxdy
\een
On the other hand,
multiply $\eta^{2\ell}w^{2q-4}w$  on both sides of (\ref{eq:NSv-2D-c}), and we have
\ben\label{eq:energy estimate-w}
II&\doteq&(2q-3)\int_{\mathbb{R}^2}|\nabla w|^2 w^{2q-4}\eta^{2\ell}dxdy\nonumber\\
&\leq&\frac{1}{2q-2}\int_{\mathbb{R}^2} w^{2q-2}\triangle (\eta^{2\ell})dxdy+\frac{1}{2q-2}\int_{\mathbb{R}^2} w^{2q-2}v\cdot \nabla(\eta^{2\ell})dxdy\nonumber\\
&&+\frac{1}{2q-2}\int_{\mathbb{R}^2} w^{2q-2}y\partial_x(\eta^{2\ell})dxdy
\een
Then it follows from (\ref{eq:w estimate-mhd4}), (\ref{eq:energy estimate-w}) that
\beno\label{eq:w estimate-mhd5}
&&\int_{\mathbb{R}^2}w^{2q}\eta^{2\ell}dxdy\\
&\leq&C(q,\ell)R^{2\alpha\beta-2}\int_{\mathbb{R}^2} w^{2q-2}\eta^{2\ell-2}dxdy+C(q,\ell)R^{3\alpha\beta-1}\int_{\mathbb{R}^2} w^{2q-2}\eta^{2\ell-1}dxdy\\
&&+C(q,\ell)R^{\alpha\beta-1}\int_{\mathbb{R}^2} |w|^{2q-1}\eta^{2\ell-1}dxdy+C(q,\ell)R^{2\alpha\beta+1-\beta}\int_{\mathbb{R}^2} w^{2q-2}\eta^{2\ell-2}dxdy\\
& \doteq& II_1+\cdots+II_4
\eeno

Noting $\ell\geq q$, by Young inequality we have
\beno
II_{1}&\leq& \delta \int_{\mathbb{R}^2}w^{2q}\eta^{(2\ell-2)\frac{q}{q-1}}dxdy+C(\delta,\ell,q)R^{1+\beta+q(2\alpha\beta-2)},
\eeno
\beno
II_{2}&\leq& \delta \int_{\mathbb{R}^2}w^{2q}\eta^{(2\ell-1)\frac{q}{q-1}}dxdy+C(\delta,\ell,q)R^{1+\beta+q(3\alpha\beta-1)},
\eeno
\beno
II_{3}&\leq& \delta \int_{\mathbb{R}^2}w^{2q}\eta^{(2\ell-1)\frac{2q}{2q-1}}dxdy+C(\delta,\ell,q)R^{1+\beta+2q(\alpha\beta-1)},
\eeno
and
\beno
II_{4}&\leq& \delta \int_{\mathbb{R}^2}w^{2q}\eta^{(2\ell)}dxdy+C(\delta,\ell,q)R^{1+\beta+q(2\alpha\beta+1-\beta)}
\eeno

Hence, firstly take
 $\delta<\frac1{32}$; secondly, for fixed $\alpha<\frac15$ and $\beta=\frac53$, we take $q_0=\frac{1+\beta}{\beta-2\alpha\beta-1}$.
Then for any $q>q_0$, we have
\beno
&&1+\beta+q(2\alpha\beta-2)<0,\quad 1+\beta+q(3\alpha\beta-1)< 0,\\
&&1+\beta+2q(\alpha\beta-1)<0, \quad 1+\beta+q(2\alpha\beta+1-\beta)< 0.
\eeno
Consequently, we get
\beno
\int_{\mathbb{R}^2}w^{2q}dxdy= 0,
\eeno
as $R\rightarrow\infty$. Thus we have $\triangle v=0$, which implies $v\equiv C$, since $v\in \chi^{\alpha,\infty}$. The proof is complete.

\section{Stability of Poiseuille flow }

In this section, we will consider the stable space of the Poiseuille flow $U = (y^2, 0)$ under the Navier-Stokes flow and prove Corollary \ref{thm:generalization of KNSS-P}.
For the Poiseuille flow $U = (y^2, 0)$, which is a solution of (\ref{eq:NS-2D}),
let $v = u-U$ be the perturbation of the
velocity, which satisfies
\begin{equation}\label{eq:NS-2D-p}
\left\{\begin{array}{llll}
-\Delta v-2+v\cdot \nabla v+\nabla \pi+(2yv_2,0)
+y^2\partial_x v=0,\\
{\rm div }~ v=0,
\end{array}\right.
\end{equation}
Let $w=\partial_2v_1-\partial_1 v_2$, then
\begin{equation}\label{eq:NSv-2D-p}
-\Delta w+v\cdot \nabla w+y^2\partial_xw=0.
\end{equation}

To overcome the singularity of the term with $y^2\partial_x$, we have to estimate the growth of the functions in $\mathcal{D}^{1,q}(\mathbb{R}^2)$.

First of all, for $1<q<2$ we have the following lemma(for example, see Theorem II.6.1 in \cite{Galdi}).
\begin{lem}\label{lem:theorem II 6 1} Let $\Omega\subseteq \mathbb{R}^2$, be an  exterior domain of locally Lipschitz and let
\beno
f\in \mathcal{D}^{1,q}(\Omega),\quad 1\leq q<2.
\eeno
Then, there exists a unique $f_0\in \mathbb{R}$ such that
\beno
f-f_0\in L^s(\Omega),\quad s=\frac{2q}{2-q}
\eeno
and for some $\gamma_1$ independent of $f$
\beno
\|f-f_0\|_s\leq \gamma_1\|\nabla(f-f_0)\|_{q}.
\eeno
\end{lem}

For the critical case $q=2$, it is obvious that $f\in BMO(\mathbb{R}^2)$ if $f\in \mathcal{D}^{1,2}(\mathbb{R}^2)$ and
$\|f\|_{*}\leq C\|u\|_{\mathcal{D}^{1,2}}$, where
\beno
\|f\|_*\doteq\sup_{x_0\in \mathbb{R}^2, r>0}\frac{1}{|Q_r(x_0)|}\int_{Q_r(x_0)}|f-f_{Q_r(x_0)}|dx<\infty,
\eeno
where $Q_r(x_0)$ is the cube whose sides have length $r$ centered at $x_0$.
It's well-known that for the BMO space, we have
\beno
\Gamma(s)=\sup_{x_0\in \mathbb{R}^2,r>0}\left(\frac{1}{|Q_r(x_0)|}\int_{Q_r(x_0)}|f-f_{Q_r(x_0)}|^sdx\right)^{\frac1s}<\infty,
\eeno
for any $s\in [1,\infty)$.

The integrable property of $f$ is stated as follows, which is a slightly different version from (1.2) in \cite{FS1972}.
\begin{lem}\label{lem:Fefferman-stein}
Let $f\in \mathcal{D}^{1,2}(\mathbb{R}^2)$. For $p\geq 1$ and $\epsilon\in (0,1)$, we have
\beno
\int_{\mathbb{R}^2}\frac{|f-f_{Q_0}|^p}{1+|x|^{2+\epsilon}}dx\leq C(\epsilon,p)\|f\|_{*}^p\leq  C(\epsilon,p)\|\nabla f\|_{L^p(\mathbb{R}^2)}^p,
\eeno
where $Q_0$ is the cube whose sides have length $1$, and is centered at the origin.

\end{lem}

{\bf Proof of Lemma \ref{lem:Fefferman-stein}.} It's similar as in in \cite{FS1972}. Let $S_k=Q_k\setminus{Q_{k-1}}$, where $Q_k$ is the cube centered at the origin whose sides have length $2^k.$
Since $\|f\|_{*}\leq C\|f\|_{\mathcal{D}^{1,2}}$, it suffices to prove that
\ben\label{eq:Ik}
I_k=\int_{S_k}\frac{|f(x)-f_{Q_0}|^p}{1+|x|^{2+\epsilon}}dx \leq C_k\|f\|_*^p,
\een
and $\sum_{k\geq 1}C_k$ is summable.

In fact,
\beno
I_k&\leq &\frac{4^{2+\epsilon}}{2^{k(2+\epsilon)}}\int_{Q_k}|f(x)-f_{Q_k}+f_{Q_k}-f_{Q_0}|^pdx \\
&\leq &  \frac{4^{2+\epsilon}}{2^{k(2+\epsilon)}} \left(2^{2k}+k^p2^{2k+2p}\right)\|f\|_*^p \leq C\frac{k^p}{2^{k\epsilon}}\|f\|_*^p
\eeno
where we used
$|f_{Q_k}-f_{Q_{k-1}}|^p\leq 4^p \|f\|_*^p. $
The proof is complete.

{\bf Proof of Corollary \ref{thm:generalization of KNSS-P}.} It's similar as the Couette flow.

(i) It's similar as Step I in the proof of Theorem \ref{thm:generalization of GW}.

For $1<q\leq 4,$
let $\eta(x,y)\in C_0^\infty(\mathbb{R}^2)$ be a cut-off function with $0\leq \eta\leq 1$ satisfying $\eta(x,y)=\eta_1(x)\eta_2(y)$, where
\begin{align*} \eta_1(x)=\left\{
\begin{aligned}
&1,\quad |x|\leq R^{\beta},\\
&0, \quad |x|>2R^{\beta}.
\end{aligned}
\right. \end{align*}
where $\beta\geq 2$, to be decided;
and
\begin{align*} \eta_2(y)=\left\{
\begin{aligned}
&1,\quad |y|\leq R,\\
&0, \quad |y|>2R.
\end{aligned}
\right. \end{align*}

For $\frac43<q_0\leq q$, multiply $q_0\eta^{2}|w|^{q_0-2}w$  on both sides of (\ref{eq:NSv-2D-p}), and we have
\ben\label{eq:energy estimate-w'}
I''&\doteq&\frac{4(q_0-1)}{q_0}\int_{\mathbb{R}^2}|\nabla (|w|^{q_0/2}|^2 \eta^{2}dxdy\nonumber\\
&\leq&\int_{\mathbb{R}^2} |w|^{q_0}\triangle (\eta^{2})dxdy+\int_{\mathbb{R}^2} |w|^{q_0}y^2\partial_x(\eta^{2})dxdy\nonumber\\
&&+\int_{\mathbb{R}^2} |w|^{q_0}v\cdot \nabla(\eta^{2})dxdy\doteq I_1''+I_2''+I_3''
\een

{\bf Case of $2<q\leq 4$. }Take $q_0=q$ and $2\leq \beta\leq (1-\frac2q)^{-1}$. Since $ v\in \mathcal{D}^{1,q}$,
obviously $I_1''\rightarrow 0$ as $R\rightarrow\infty.$
 When $\beta\geq 2$,
\beno
I_2''\leq C(q)R^{2-\beta}\int_{R^\beta\leq |x|\leq 2R^\beta}|w|^qdxdy\rightarrow 0
\eeno
About the term $I_3$, due to Lemma \ref{lem:Galdi}, we have
\beno
|v(x,y)|\leq |(x,y)|^{1-\frac{2}{q}}
\eeno
Then we have
\beno
I_3''\leq C(q)R^{\beta(1-\frac2q)-1}\int_{R^\beta\leq |x|\leq 2R^\beta}|w|^qdxdy\rightarrow 0,
\eeno
as $R\rightarrow\infty$, since
\beno
\beta(1-\frac2q)\leq 1.
\eeno

Consequently, we get
$
\nabla (|w|^{\frac{q}{2}})\equiv0
$
which implies that
$
w\equiv C.
$
The same arguments hold.

{\bf Case of $q=2$. } Take $q_0=\frac95$ and $\beta=3$. Obviously $I_1''+I_2''\rightarrow 0$ as $R\rightarrow\infty$ by H\"{o}lder inequality. Next, we estimate the term $I_3''$. With the help of Lemma \ref{lem:Fefferman-stein} as $\epsilon=1$,
there holds
\beno
\int_{\mathbb{R}^2}\frac{|v-v_{Q_0}|^{10}}{1+|(x,y)|^{3}}dxdy<C\|v\|_{*}^{10}
\eeno
Thus
\beno
I_3''\leq C(q)R^{-1+\frac{3}{10}\beta}\rightarrow 0,
\eeno
as $R\rightarrow\infty$.

{\bf Case of $\frac43<q<2$. } At this time, take $q_0=\frac{3q-2}{2}\in (1,q)$ and $\beta=1+\frac{2}{q_0}\in (2,3)$. Then
\beno
2-\beta+(1+\beta)(1-\frac{q_0}{q} )=0
\eeno
and hence
\beno
I_2''\leq C(q)\left(\int_{R^\beta\leq |x|\leq 2R^\beta}|w|^qdxdy\right)^{\frac{q_0}{q}}R^{2-\beta+(1+\beta)(1-\frac{q_0}{q} )}\rightarrow0
\eeno
as $R\rightarrow\infty.$

For the term $I_3''$,
by Lemma \ref{lem:theorem II 6 1} there exists  a constant vector $v_0$ such that
\beno
\|v-v_0\|_{L^{\frac{2q}{2-q}}}\leq C\|\nabla v\|_{L^q}.
\eeno
Thus
\beno
I_3''\leq C(q)\left(\int_{R^\beta\leq |x|\leq 2R^\beta}|w|^qdxdy\right)^{\frac{q_0}{q}}\|\nabla v\|_{L^q} R^{-1+(1+\beta)(1-\frac{q_0}{q}-\frac{2-q}{2q} )}\rightarrow0
\eeno
as $R\rightarrow\infty,$ since
\beno
1-\frac{q_0}{q}-\frac{2-q}{2q}=0.
\eeno

(ii)
Then it follows from (\ref{eq:w estimate-mhd5}) that
\beno\label{eq:w estimate-pois}
&&\int_{\mathbb{R}^2}w^{2q}\eta^{2\ell}dx\\
&\leq&C(q,\ell)R^{2\alpha\beta-2}\int_{\mathbb{R}^2} w^{2q-2}(\eta^{2\ell-2})dxdy+C(q,\ell)R^{3\alpha\beta-1}\int_{\mathbb{R}^2} w^{2q-2}(\eta^{2\ell-1})dxdy\\
&&+C(q,\ell)R^{\alpha\beta-1}\int_{\mathbb{R}^2} |w|^{2q-1}(\eta^{2\ell-1})dxdy+C(q,\ell)R^{2\alpha\beta+2-\beta}\int_{\mathbb{R}^2} w^{2q-2}(\eta^{2\ell-2})dxdy\\
&&  =II_1'+\cdots+II_4'
\eeno
At this time, since $0<\alpha<\frac18$, we can
choose
\beno
\beta=\frac83,
\eeno
thus
\beno
1+\beta+q(3\alpha\beta-1)< 0,1+\beta+q(2\alpha\beta+2-\beta)< 0,
\eeno
for $q$ sufficiently large.
Since the similar arguments as Theorem \ref{thm:generalization of KNSS-C} hold,
we omitted it. The proof is complete.

\noindent {\bf Acknowledgments.}
W. Wang was supported by NSFC under grant 11671067 and
 "the Fundamental Research Funds for the Central Universities".

\end{document}